\documentclass[10pt]{amsart}
\usepackage{amscd} 
\usepackage{amsfonts} 
\usepackage{amssymb} 
\usepackage{latexsym}

\newcommand{\ncm}{\newcommand}

\newtheorem{thm}{Theorem}[section]
\newtheorem{prop}[thm]{Proposition}
\newtheorem{lemma}[thm]{Lemma}
\newtheorem{cor}[thm]{Corollary}
\newtheorem{lem&def}[thm]{Lemma \& Definition}
\newtheorem{defi}[thm]{Definition}

\def\C{\mathbb{C}\,}

\def\H{\mathbb{H}\,}

\ncm{\End}{\mbox{\rm End}\,}
\def\|{\, | \, }

\def\Hom{\mbox{\rm Hom}\,}

\def\id{\mbox{\rm id}}

\def\Im{\mbox{\rm Im}\,}

\def\to{\rightarrow}

\def\o{\otimes}    
\def\bra{\langle}
\def\ket{\rangle}

\ncm{\rarr}[1]{\stackrel{#1}{\longrightarrow}}
\ncm{\larr}[1]{\stackrel{#1}{\longleftarrow}}
\def\i{^{(1)}}
\def\j{^{(2)}}

\def\cop{\Delta}

\def\eps{\varepsilon}

\def\du1{\hat 1}

\def\0{_{(0)}}
\def\1{_{(1)}}
\def\2{_{(2)}}
\def\3{_{(3)}}
\def\t{^{\tau}}

\def\du1{\hat 1}
\def\lact{\triangleright}
\def\ract{\triangleleft}

\begin{document}

\title[Hopf Algebroids and Irreducible Extensions]{Skew Hopf algebras, Irreducible Extensions
and the $\Pi$-Method}
\author{Lars Kadison} 
\address{Department of Mathematics \\ University of Pennsylvania \\
David Rittenhouse Lab, 209 S. 33rd St. \\ 
Philadelphia, PA 19104 \\
currently: Louisiana State University,
Baton Rouge, LA 70803} 
\email{lkadison@math.upenn.edu}
\thanks{}
\date{} 

\begin{abstract}
To a depth two extension $A | B$, we associate the 
dual
bialgebroids $S := \End\,_BA_B$ and $T := (A \o_B A)^B$ over the centralizer 
$R=C_A(B)$.  In a set-up which is quite
common, where $R$ is a 
subalgebra of $B$,  two nondegenerate pairings of $S$ and $T$ will 
define an anti-automorphism $\tau$ of the algebra $S$.  
Making use of a two-sided depth two structure, we
 show that $\tau$ is an antipode and $S$ is  a Hopf algebroid of a type we
call skew Hopf algebra.  A final section discusses how $\tau$ and the nondegenerate pairings generalize to modules via the $\pi$-method for depth two.
\end{abstract} 
\maketitle

\section{Introduction}

For reasons of symmetry in representation theory, given a bialgebra or bialgebroid one would like to expose the presence
of an antipode. In the reconstruction of
Hopf algebras and weak Hopf algebras in subfactor theory a key idea in the definition of antipode is  to make use of the existence of two nondegenerate pairings of dual bialgebras. In terms of
a depth two, finite index subfactor $N \subseteq M$ with trivial relative commutant, 
its basic construction $M \subseteq M_1$,
and another one above, $M_1 \subseteq M_2$, there are
conditional expectations $E_M: M_1 \rightarrow M$ and $E_{M_1}: M_2 
\rightarrow M_1$.  In addition to  the Jones projections
$e_1 \in M_1$ and $e_2 \in M_2$, the two relative commutants that are paired
nondegenerately are
in ordinary algebraic centralizer notation $C = C_{M_1}(N)$
and $V = C_{M_2}(M)$.  The antipode $\tau: V \rightarrow V$ is then
defined as the ``difference'' of two such pairings:
\begin{equation}
E_M E_{M_1}(ve_1 e_2 c) = E_M E_{M_1}(ce_2 e_1 \tau(v))
\end{equation}
for $c \in C$ and $v \in V$:  see for example, \cite{KN,  NV1, S}
for the details of why this formula works.  

 The duality method for defining antipode has been lying dormant in recent generalizations
of depth two to algebras and rings and actions of bialgebroids on these. For example,  antipodes have been defined recently 
in the case of a Frobenius extension $A | B$ as the restriction of a standard
anti-isomorphism of     the left and right endomorphism rings
$\End A_B \rightarrow \End {}_BA$ to an anti-automorphism of the subring of bimodule endomorphisms
$S = \End {}_BA_B$:  on depth two Frobenius extension this defines the antipode or its inverse in \cite{BS}
in a  dual way on both $S$ and $T$; it also necessitates a revision of the definition of the notion of Hopf algebroid
using the notions of left and right bialgebroid. 
 Antipodes have also been defined from geometric ideas of Lu \cite{Lu} 
for H-separable extensions \cite{LK}, extensions of Kanzaki 
separable algebras \cite{LK2006A} and Hopf-Galois extensions \cite{LK2005},
and from group theory in \cite{LK2007} for pseudo-Galois extensions.

In this paper we define antipode as the difference of two nondegenerate
pairings for a special extension $A | B$ where the centralizer $R$ is a 
subalgebra of the smaller algebra $B$, which we call irreducible extension.  
In this case, the two hom-groups, the left and right $R$-duals of the
bialgebroid $T$, coincide. The plan to find an antipode from
this identity and 
satisfying the several axioms of a Hopf algebroid works
 not so much because of any Frobenius structure (as assumed previously)
but on a two-sided depth two structure as shown in sections~3 and~4 below.  
We recall that the dual bialgebroids $S$ and $T$ depend only a
one-sided depth two structure \cite{LK2006A}, but at two stages in this
paper (using both nondegenerate pairing and the proposition in section 3) we require a two-sided depth two
structure.  The Frobenius extension hypothesis avoided in this paper,
makes one-sided depth two extensions two-sided \cite[6.4]{KS}.
The Hopf algebroid structure we obtain on $S$ is very nearly a
Hopf algebra which is finite projective over a commutative base ring:
we discuss its properties after Theorem~\ref{thm-ha} and designate as skew
Hopf algebras such a
Hopf algebroid.
We end with
a discussion of how $\tau$ and the nondegenerate pairings generalize to modules via the $\pi$-method for depth two.
A certain mapping between cochain complexes
formed from the left- and right-handed $\pi$-methods is shown to be nullhomotopic.

\section{Preliminaries on depth two extensions}

Let $B$ be a unital subalgebra of $A$, an associative noncommutative algebra with unit
over a commutative ground ring $K$.  The algebra extension $A | B$ is \textit{depth two} if there is
a positive integer $N$ such that 
\begin{equation}
\label{eq: D2}
A \o_B A \oplus * \cong  A^N
\end{equation}
 as natural $B$-$A$ (left D2) and $A$-$B$-bimodules (right D2)
\cite{KS}. For an $A$-$A$-bimodule $M$, denote the subgroup of $B$-central elements in $M$
by $$M^B := \{ m \in M | \forall b \in B, bm = mb \}.$$  
Equivalently, the algebra extension $A | B$ is depth two if there are elements 
$$\beta_i \in S := \End\,_BA_B, \ \ t_i \in T := (A \o_B A)^B$$
 (called a \textit{left D2 quasibasis}) such that all elements of 
$A \o_B A$ may be written as ($a,a' \in A$)
\begin{equation}
\label{eq: left d2 qb}
 a \o a' = \sum_{i=1}^N t_i\beta_i(a)a', 
\end{equation}
and similarly, there are elements (of a right D2 quasibasis) $\gamma_j \in S$, $u_j \in T$ such that
\begin{equation}
\label{eq: right d2 qb}
a \o a' = \sum_j a \gamma_j(a') u_j. 
\end{equation}
For example, given a right D2 quasibasis, define
a split $A$-$B$-epimorphism $A^N \rightarrow A \o_B A$ by
$$ (a_1, \ldots, a_N) \longmapsto \sum_{j = 1}^N a_j u_j $$
which is split by the $A$-$B$-monomorphism $A \o_B A \rightarrow A^N$ given
by $x \o_B y \mapsto (x\gamma_1(y), \ldots, x \gamma_N(y))$.  
Conversely, given a split epi $A^N \rightarrow A \o_B A$, we obtain
mappings $\sum_{i=1}^N f_i \circ g_i$ where $g_i \in \Hom (A \o_B A, A)$
and $f_i \in \Hom (A, A \o_B A)$; but there are somewhat obvious isomorphisms
$\Hom (A, A \o_B A) \cong T$ and $\Hom (A \o_B A, A) \cong S$ in either 
case of $A$-$B$- or $B$-$A$-bimodule homomorphisms.  
 We fix the notations for both right and left D2 quasibases
throughout this paper.

For example, an H-separable extension $A | B$ is of depth two since
the condition above on the tensor-square holds even more strongly as natural 
$A$-$A$-bimodules. Another example:  $A$ a f.g.\ projective
algebra over commutative ground ring $B$, since  left or right D2 quasibases
are easily constructed from a dual basis.
As a third class of  examples, consider a Hopf-Galois extension $A | B$ with
$n$-dimensional Hopf $k$-algebra $H$ \cite{KT}.  Recall that $H$ acts from the 
left on $A$ with subalgebra of
invariants $B$, induces
a dual right coaction $A \to A \o_k H^*$, $a \mapsto a\0 \o a\1$,
 and Galois isomorphism $\beta: A \o_B A \stackrel{\cong}{\longrightarrow}
A \o_k H^*$ given by $\beta(a \o a') = a{a'}\0 \o {a'}\1$ , 
which is an $A$-$B$-bimodule, right $H^*$-comodule morphism.
It follows that $A \o_B A \cong \oplus^n A$ as $A$-$B$-bimodules;
as $B$-$A$-bimodules there is a similar isomorphism
 by making use of the alternative
 Galois isomorphism $\beta'$ given by  $\beta'(a \o a') = a\0 a' \o a\1$. 
The paper \cite{LK2006} extends the definition above of depth two to include
the case where the tensor-square of $A | B$ is isomorphic to any direct sum of
$A$ with itself (not necessarily a finite direct sum as in eq.~\ref{eq: D2}); thus any Hopf-Galois extension
is depth two in this extended sense.  However, this theory does not have a theory
of dual bialgebroids congenial for the results in this paper, and we shall not make
use of it.    

The papers  \cite{KS, LK} defined dual bialgebroids with action and smash product structure 
within the endomorphism ring tower construction above a depth two ring extension $A | B$.  In more
detail, if $R$ denotes the centralizer of $B$ in $A$,
a left $R$-bialgebroid structure on $S$ is given by the composition
ring structure on $S$ with source and target mappings
corresponding to the left regular representation $\lambda: R \to S$
and right regular representation $\rho: R^{\rm op} \to S$, respectively. Since these
 commute
($\lambda_r \rho_s = \rho_s \lambda_r$ for every $r,s \in R$),
we may induce an $R$-bimodule structure on $S$ solely
  from the left  
by $$r \cdot \alpha \cdot s := \lambda_r\rho_s \alpha = r\alpha(-)s.$$

Now an $R$-coring structure $(S, \cop, \eps)$ is given
by 
\begin{equation}
\label{eq: cop A}
\cop(\alpha) := \sum_i \alpha(-t_i^1) t_i^2 \o_R \beta_i 
\end{equation}
for every $\alpha \in S$, denoting $t_i = t_i^1 \o_B t_i^2 \in T$ by
suppressing a possible summation, and
\begin{equation}
\label{eq: eps A}
\eps(\alpha) = \alpha(1)
\end{equation}  
satisfying the additional axioms of a bialgebroid (cf.\ appendix),
such as multiplicativity of $\cop$ and a condition that makes sense of this
requirement. We have the equivalent formula for the coproduct \cite[Th'm 4.1]{KS}:
\begin{equation}
\label{eq: right D2 cop}
\cop(\alpha) := \sum_j \gamma_j \o_R u_j^1 \alpha(u_j^2 -) 
\end{equation}

For a depth two extension, 
\begin{equation}
\label{eq: id}
S \o_R S \stackrel{\cong}{\longrightarrow} \Hom ({}_BA \o_B A_B, {}_BA_B), \ \
\alpha \o_R \beta \longmapsto (x \o_B y \mapsto \alpha(x)\beta(y))
\end{equation}
with inverse provided by $F \mapsto \sum_i F(-\o t_i^1)t_i^2 \o_R \beta_i$.   
Under this identification, the formula for the coproduct $\cop: S \rightarrow S \o_R S$
becomes
\begin{equation}
\alpha\1(x) \alpha\2(y) = \alpha(xy)
\end{equation}
using either equation.  

The left action of $S$ on $A$ given by evaluation, $\alpha \lact a = \alpha(a)$,
has invariant subalgebra (of elements
$a \in A$ such that $\alpha\lact a = \eps(\alpha)a$) equal precisely to $B$ if the natural module $A_B$ is balanced \cite{KS}. 
This action is measuring  by
eq.~(\ref{eq: right D2 cop}). 

The smash product $A \rtimes S$, which is $A \o_R S$ as  abelian groups with associative
multiplication given by
\begin{equation}
(x \rtimes \alpha)(y \rtimes \beta) = x (\alpha\1 \lact y) \rtimes \alpha\2 \beta,
\end{equation}
is isomorphic
as rings to $\End A_B$ via $a \o_R \alpha \mapsto \lambda_a \alpha$ \cite{KS}.

In general $T = (A \o_B A)^B$ has a unital ring structure induced from 
$T \cong \End\,_A(A\! \o_B\! A)_A  $ via $F \mapsto F(1 \o 1)$, which is given
by 
\begin{equation}
 tu = {u}^1 t^1 \o t^2 {u}^2
\end{equation} for each $t,u \in T$. There are
obvious commuting homomorphisms of $R$ and $R^{\rm op}$ into $T$ given by 
$r \mapsto 1 \o r$ and $s \mapsto s \o 1$, respectively. From the right,
these two source and target mappings induce the $R$-$R$-bimodule
structure ${}_RT_R$ given by $$r \cdot t \cdot s =  (t^1 \o t^2)(r \o s) =
rt^1 \o t^2s,$$
the ordinary bimodule structure on a tensor product.  

There is a right $R$-bialgebroid structure on $T$ with coring structure
$(T, \cop, \eps)$ given by the two equivalent formulas: 
\begin{equation}
\label{eq: copB}
\cop(t) = \sum_i t_i \o_R (\beta_i(t^1) \o_B t^2) = \sum_j (t^1 \o_B \gamma_j(t^2)) \o_R
u_j 
\end{equation}
\begin{equation}
\label{eq: epsB}
\eps(t) = t^1 t^2
\end{equation}
By \cite[Th'm 5.2]{KS} $\cop$ is multiplicative and the other axioms of
a right bialgebroid are satisfied. 

As an example of $S$ and $T$, consider the Hopf-Galois extension $A | B$ of
$k$-algebras introduced
above.  Since $\beta$ is an $A$-$B$-isomorphism, we may compute
that $T \cong R \o_k H^*$
via $\beta$, which induces a smash product structure on $R \o H^*$ relative to the
Miyashta-Ulbrich action of $H^*$ on $R$ from the right. 
The well-known isomorphism $\End A_B \cong A \rtimes H$ via $a \rtimes h \mapsto
\lambda(a)(h \lact \cdot)$ restricts to $S \cong R \rtimes H$, i.e., $S$
is a smash product of $R$ with $H$ via the restriction of the left action of $H$
to $R$.  In both cases, the $R$-coring structures are the trivial ones induced
from the coalgebras $H$ and $H^*$.  

There is a right action of $T$ on $\mathcal{E} := \End\,_BA$ 
given by $f \ract t = t^1 f(t^2 -)$ for $f \in \mathcal{E}$. This is
a measuring action by Eq.\ (\ref{eq: left d2 qb}) since
$$ (f \ract t\1)\circ (g \ract t\2) = \sum_i t^1_i f(t^2_i\beta_i(t^1)g(t^2 -)) = fg \ract t. $$
The subring of invariants in $\mathcal{E}$ is $\rho(A)$ \cite{KS}.
Also, 
in analogy with $\End A_B \cong A \rtimes S$, the smash product
ring $T \ltimes \mathcal{E}$ is isomorphic to $\End\,_AA\! \o_B\! A$ via
$\Psi$ given by 
\begin{equation}
\Psi(t \o f)( a \o a') = at^1 \o_B t^2 f(a').
\end{equation}

Sweedler \cite{Sw} defines left and right $R$-dual rings of an $R$-coring.  In
the case of a left $R$-bialgebroid $H$ with $H_R$ and ${}_RH$
finitely generated projective, such as $(S,\lambda, \rho, \cop, \eps)$ above,
the left and right Sweedler $R$-dual rings are extended to right bialgebroids $H^*$
and ${}^*\!H$ in \cite{KS}. For example, $H^*$ has a natural nondegenerate pairing
with $H$ denoted by $\bra h^*, h \ket \in R$ for $h^* \in H^*, h \in H$. Then 
the $R$-bimodule structure on $H^*$, multiplication,
and comultiplication are given below, respectively, where
$R \stackrel{s}{\to} H \stackrel{t}{\leftarrow} R^{\rm op}$ denotes
the commuting morphism set-up of the bialgebroid $H$:
\begin{eqnarray}
\label{eq: dual 1}
\bra r \cdot h^* \cdot r', h \ket & := & r \bra h^*, ht(r') \ket \\
\bra h^* g^*, h \ket & := & \bra g^*, \bra h^*, h\1 \ket \cdot h\2 \ket 
\label{eq: Sweedlers} \\ \label{eq: dual 2}
\bra h^*, hh' \ket & := & \bra {h^*}\1 \cdot \bra {h^*}\2, h' \ket, h \ket 
\end{eqnarray}
Of course, the unit of $H^*$ is
$\eps_H$ while the counit on $H^*$ is  $\eps(h^*) = \bra h^*, 1_H \ket$. 
Eq.\ (\ref{eq: Sweedlers}) is the formula for multiplication \cite[3.2(b)]{Sw}.

There are similar formulas for the right bialgebroid structure on the left
$R$-dual ${}^*\!H$: see \cite[2.6]{KS}.  In the particular case of the left bialgebroid
$S$ of a depth two ring extension, it turns out that $S$ is isomorphic as $R$-bialgebroids
 to both $R$-duals, $T^*$ and ${}^*\!T$ via two  nondegenerate pairings,
one of which is given by 
($\alpha \in S, t \in T$): 
\begin{equation}
\label{eq: bra-ket pairing}
\bra  \alpha, t \ket =  \alpha (t^1)t^2 \in R
\end{equation}
This induces an isomorphism of left $R$-modules $S \rightarrow \Hom (T_R, R_R)$
via $\alpha \mapsto \bra \alpha, - \ket$ with inverse $$\phi \mapsto \sum_i \phi(t_i)\beta_i.$$
Significantly, there is another
nondegenerate pairing of $S$ and $T$ for left and right D2 extensions given by
\begin{equation}
\label{eq: bracket pairing}
[ t, \alpha ] = t^1 \alpha( t^2)
\end{equation}
  \cite[5.3]{KS}.   This induces an isomorphism of right $R$-modules $S \rightarrow
\Hom ({}_RT, {}_RR)$ given by $\alpha \mapsto [ - , \alpha]$, with inverse
given from a right D2 quasibasis by $$\psi \mapsto \sum_j \gamma_j(-) \psi(u_j).$$


\section{Irreducible extensions}
\label{sec: irred}

We define a class of depth two extension where we may readily exploit 
the two nondegenerate pairings just given in eqs.~(\ref{eq: bra-ket pairing})
and~(\ref{eq: bracket pairing}).  We say that an algebra extension $A | B$ is \textit{irreducible} if it is depth two and its centralizer $C_A(B) = R$ is a subalgebra
of $B$, so $R \subseteq B$. Then $R$ is a commutative subalgebra,
since $r s = sr$ for all $s,r \in R$ follows from noting for instance
that $r \in B$ and $s \in A^B$. 

 This set-up is
quite common. For example, irreducible depth two subfactors are irreducible in our sense since the centralizer
is one-dimensional over the complex numbers \cite{KN, NV1}. A second example: Taft's Hopf algebras including Sweedler's $4$-dimensional Hopf algebra, which are 
generated by a grouplike element $g$ and a skew-primitive element $x$, over the commutative Frobenius
subalgebra $B$ generated by the element $x$:
this extension satisfies $B = R$ and is depth two
(in fact strongly graded, therefore Hopf-Galois) \cite{GGG, K}. A third example
is the extension $\C \subset \H$, the complex numbers as a subring in the
real quaternions.   

Another type of
example of irreducible extension is the
 H-separable extension of full $n \times n$ matrix algebra over the triangular matrix subalgebra, which of course 
has trivial centralizer. (If we pass to  infinite dimensional matrices of finite type, a version of this example shows
that Cuadra's result for separable Hopf-Galois extension \cite{C} does not extend
to separable, depth two extensions; namely, an example of a non-finitely generated H-separable
extension.)
Finally, note that an intermediate ring
$B$ in an 
irreducible extension $A \| C$  is irreducible if $A \| B$ is D2, since
$A^B \subseteq A^C \subseteq C \subseteq B$. 

   For an irreducible extension $A | B$ with the construction $T = (A \o_B A)^B$, 
we note that
\begin{equation}
\Hom (T_R, R_R) = \Hom ({}_RT, {}_RR).
\end{equation}
This follows from $R \subseteq B$ and commutativity in $R$, for given $\phi \in \Hom (T_R, R_R),
t \in T, r \in R$:
$$\phi(rt) = \phi(tr) = \phi(t)r = r\phi(t),$$
and a similar computation showing ${}^*T \subseteq T^*$ using left and right $R$-dual notation.
In case $A | B$ is (two-sided) depth two and irreducible, the two nondegenerate pairings
($\alpha, \beta \in S$) 
$$S \stackrel{\cong}{\longrightarrow} \Hom (T_R, R_R), \ \ \alpha \longmapsto \bra \alpha, - \ket$$
and $$ S \stackrel{\cong}{\longrightarrow} \Hom ({}_RT, {}_RR), \ \  \beta \longmapsto [-, \beta]$$
induce a bijection $\tau$ of $S$ with itself completing a commutative
 triangle with these two mappings.
Then define 
\begin{equation}
\label{eq: tau}
\tau: S \rightarrow S, \ \ \ \bra \alpha, t \ket = [t, \tau(\alpha) ] 
\end{equation}
for all $t \in T$ and $\alpha \in S$.  We also make use of the notation $\alpha^{\tau} = \tau(\alpha)$,
for which the last equation becomes
\begin{equation}
\label{eq: tau alt}
\alpha(t^1)t^2 = t^1 \alpha^{\tau}(t^2).
\end{equation}

Notice that this approach will not work on the two nondegenerate
pairings $T \rightarrow \Hom (S,R)$, to define a self-bijection
on $T$, unless we assume that $R$
coincides with the center of $A$.

\begin{lemma}
The mapping $\tau: S \rightarrow S$ is an anti-automorphism of $S$
satisfying $\tau(\rho_r) = \lambda_r$, $\tau(\lambda_r) = \rho_r$
for $r \in R$ and $\alpha\t (1) = \alpha(1) $.
\end{lemma}
\begin{proof}
It is clear that $\tau$ is linear and bijective.  
We note that for $t = t^1 \o_B t^2 \in T$, $\alpha, \beta \in S$, 
$$ t^1 \beta^{\tau} \alpha^{\tau}(t^2) = \beta(t^1) \alpha^{\tau}(t^2) = \alpha\beta(t^1)t^2 =
t^1 (\alpha \beta)^{\tau}(t^2) $$
since $t^1 \o_B \alpha^{\tau}(t^2)$ and $\beta(t^1) \o_B t^2$ both are in $T$.  
Then $[t, \beta^{\tau} \alpha^{\tau}] = [ t, (\alpha \beta)^{\tau}]$,
so by nondegeneracy of this pairing, $\tau$ is an anti-automorphism of $S$. 

We also check that
$$ [ t, \rho_r^{\tau}] = \bra \rho_r, t \ket = t^1 r t^2 = [t, \lambda_r] $$
whence $\rho_r^{\tau} = \lambda_r$.  Since $R \subset B$ and $T = (A \o_B A)^B$,
we note that
$$ [t, \rho_r] = t^1t^2 r = rt^1 t^2 = \bra \lambda_r, t \ket = [t, \lambda_r^{\tau}] $$
whence $\lambda_r^{\tau} = \rho_r$ for each $r \in R$.

Finally, with $1_T = 1 \o_B 1 $, the equality $\alpha(1) = \alpha\t(1)$
follows from $\bra \alpha, 1_T \ket = [ 1_T, \alpha\t ]$.   
\end{proof}

\begin{prop}
Suppose $A | B$ is an irreducible extension
with left D2 quasibasis $t_i \in T, \beta_i \in S$ and right D2 quasibasis
$u_j \in T, \gamma_j \in S$. Then $t_i, \beta_i^{\tau}$ is a right D2 quasibasis
and $u_j, \gamma_j^{\tau}$ is a left D2 quasibasis.
\end{prop}
\begin{proof}
Note that for $t \in T$, a special instance of eq.~\ref{eq: left d2 qb} yields
$$t = \sum_i t_i^1 \o_B t^2_i \beta_i(t^1)t^2 = \sum_i t^1 \beta_i^{\tau}(t^2) t_i^1 \o_B t_i^2,$$
since $\beta_i(t^1)t^2 \in R \subseteq B$. Now recall that
$ A \o_R T \cong A \o_B A$ via $a \o_R t \mapsto at^1 \o_B t^2$
since $A | B$ is right D2 (for an inverse is given by $x \o_B y \mapsto
\sum_j x \gamma_j(y) \o_R u_j$).  Note that
$$ A \o_B A \longrightarrow A \o_R T, \ \ x \o_B y \longmapsto \sum_i x \beta_i^{\tau}(y) \o_R t_i$$
is a left inverse of $A \o_R T \rightarrow A \o_B A$, $a \o_R t \mapsto at^1 \o_B t^2$
since
$$ \sum_i at^1 \beta_i^{\tau}(t^2) \o_R t_i = a \o_R \sum_i t^1 \beta_i^{\tau}(t^2)t_i = a \o_R t.$$
Hence, it is also a right inverse, so
\begin{equation}
x \o_B y = \sum_i x \beta_i^{\tau}(y)t_i
\end{equation}  
for all $x, y \in A$, which shows that $t_i \in T, \beta_i^{\tau} \in S$ is
a right D2 quasibasis.  

The argument that $u_j, \gamma_j^{\tau}$ is a left D2 quasibasis is very similar.  
\end{proof}


\section{The Hopf algebroid $S$}

We are now in a position to show that the anti-automorphism $\tau$ on
$S$,
defined in eq.~(\ref{eq: tau alt}), is an antipode satisfying the axioms of B\"ohm-Szlach\'anyi \cite{BS}.
In order for $S$ to be a Hopf algebroid in the sense of Lu, we need
one additional requirement, e.g. that $R$ be a separable $K$-algebra,
in order that we may find a section of the canonical epi $S \o_K S \rightarrow S \o_R S$.  

\begin{thm}
\label{thm-ha}
If $A | B$ is an irreducible  extension, then the anti-automorphism
$\tau$ on $S = \End {}_BA_B$
is an antipode and $S$ is a Hopf algebroid.
\end{thm}
\begin{proof}
We check that the axioms of a Hopf algebroid are satisfied,
axioms given in for example \cite[8.7]{KS} and repeated in an appendix below.
Define a right bialgebroid structure on $S$ over $R^{\rm op} = R$ by
choosing target map $t_R = \lambda$ and source map $s_R = \rho$, whence
the $R$-bimodule structure on $S$ becomes
$$ r \cdot \alpha \cdot s = \alpha \rho_s \lambda_r = \alpha(r ? s) = r \alpha(?)s,$$
the usual structure on $S$ introduced above (since $R \subseteq B$).  
Then the left bialgebroid structure $(S, R, \cop, \eps)$ introduced above
in eqs.~(\ref{eq: cop A}),~(\ref{eq: eps A}) and~(\ref{eq: right D2 cop}) is also a right bialgebroid
structure.  In other words, the axioms (1) and (2) in \cite[8.7]{KS} are satisfied
by noting that $s_L = t_R$,
$t_L = s_R$, and taking $\cop_L = \cop_R$ and $\eps_L = \eps_R$. 
We check that also the axiom of a right bialgebroid,
$s_R(r)\alpha\1 \o_R \alpha\2 = \alpha\1 \o_R t_R(r) \alpha\2 $
is satisfied since
$\rho_r \alpha\1(x)\alpha\2(y) = \alpha\1(x) \lambda_r \alpha\2(y)$ ($x,y \in A$)  
in the identification $S \o_R S \cong \Hom ({}_BA \o_B A_B, {}_BA_B)$ 
in eq.~(\ref{eq: id}).  In addition, the axiom
$\eps(t_R(\eps(\alpha))\beta) = \eps(\alpha \beta) = \eps(s_R(\eps(\alpha)\beta)$
($\alpha, \beta \in S$) is satisfied since
$$\eps(\lambda_{\eps(\alpha)}\beta) = \alpha(1)\beta(1) = \alpha(\beta(1)) = \eps(\alpha \beta)$$
which equals $\beta(1)\alpha(1) = \eps(\rho_{\eps(\alpha)}\beta)$,
where we use $\alpha(1) \in R \subseteq B$ and $R$ is commutative.

We proceed to axiom (i):
$$ \tau(\alpha t_R(r)) = \tau(\lambda_r)\tau(\alpha) = s_R(r) \tau(\alpha)$$
by the lemma, and
$$ \tau(t_L(r)\alpha) = \tau(\alpha)\tau(\rho_r) = \tau(\alpha)s_L(r) $$
for all $r \in R, \alpha \in S$.  

Finally, axiom (ii) is satisfied since by eq.~(\ref{eq: cop A})
$$ \alpha\1 \circ \tau(\alpha\2) = \sum_i \alpha(\beta_i^{\tau}(-) t_i^1)t_i^2 = \lambda_{\alpha(1)}
  = s_L(\eps(\alpha))$$
by the proposition, and by lemma, 
$$\tau(\alpha\1)\circ \alpha\2 = 
\sum_i \tau(\rho_{t^2_i} \alpha \rho_{t^1_i}) \circ \beta_i =
\sum_i \lambda_{t^1_i} \alpha\t \lambda_{t^2_i} \beta_i $$
$$ = \rho_{\alpha\t(1)} = s_R(\eps(\alpha)).$$
Thus $S$ and $\tau$ form a Hopf algebroid. 
\end{proof}

For example, the Hopf algebroid structure on $S$ coincides with
that in \cite{LK} should the extension be H-separable as well 
as irreducible, since $\tau$ exchanges
$\lambda_r$ and $\rho_r$.

Note that Hopf algebroid $(S,\tau)$ satisfies properties
close to a Hopf algebra, among them:
\begin{enumerate}
\item $\tau(r \cdot \alpha) = \tau(\alpha) \cdot r$ and $\tau(\alpha \cdot r) = r \cdot \tau(\alpha)$ for all $r \in R, \alpha \in S$; 
\item $\eps(\alpha \beta) = \eps(\alpha)\eps(\beta)$ for all $\alpha, \beta \in S$;
\item $\alpha\1 \tau(\alpha\2) = \eps(\alpha) \cdot 1_S$ for all $\alpha \in S$;
\item $\tau(\alpha\1)\alpha\2 = 1_S \cdot \eps(\alpha)$ for all $\alpha \in S$;
\item $\eps \tau = \eps$;
\item $\tau$ is an ``anti-coalgebra homomorphism.''
\end{enumerate}
Also, $S$ is finite projective over the base ring $R$, which is commutative.
However, such basic algebraic properties of a Hopf algebra as
$r \cdot 1 = 1 \cdot r$ and $(\alpha \cdot r) \beta = \alpha (r \cdot \beta)$
are suspended for $S$ (unless $R$ coincides with the center of $A$).
We propose to call a Hopf algebroid with equal right and left bialgebroid structures over a commutative base ring,
possessing an anti-automorphism
exchanging source and target,
both mappings with image in the center,  and satisfying the properties enumerated directly above, a \textit{skew Hopf algebra}.  Given the general
nature of the example $S$, we  would expect that skew Hopf algebras are quite
common occurences.   
\begin{cor}
Suppose $A | B$ is an irreducible extension and
$R$ is a separable $K$-algebra.  Then ($S$, $\tau$) is a Hopf algebroid
in the sense of Lu.
\end{cor}
\begin{proof}
Let $e = e^1 \o_K e^2$ be a separability element for $R$.
Define a section $\eta: S \o_R S \rightarrow S \o_K S$
of the canonical epi $S \o_K S \rightarrow S \o_R S$
by $\eta(\alpha \o_R \beta) = \alpha \cdot e^1 \o_K e^2 \cdot \beta$,
since $e^1 e^2 = 1$ and $re = er$ for $r \in R$.  
Then the axiom $\mu (\id \o \tau) \eta \cop = s_L \circ \eps$
follows from
$$[\mu (\id \o \tau) \eta \cop(\alpha), u ] = 
\sum_i u^1 \alpha({\beta_i}\t(u^2e^2)t_i^1)t^2_i e^1
= u^1 \alpha(1) u^2 = [\lambda_{\alpha(1)}, u]$$
by the proposition and since $R \subseteq B$, $e^2 e^1 = 1$. 
The other axioms follow as in the proof of the theorem. 
\end{proof}
   
Of course, like the inverse in group theory, antipodes are important to
the representation theory of a bialgebroid.  For example, we can now 
define a right $S$-module algebra structure on $A$ from the left structure
 by $a \ract \alpha =
\alpha\t \lact a$, which satisfies the measuring rule $(xy) \ract \alpha
= (x \ract \alpha\2)(y \ract \alpha\1)$ for $x,y \in A$ and $\alpha \in S$.    

The antipode on $S$ will not dualize readily to an antipode on $T$ without
the duality properties discussed in \cite{B}, such as $S$ possessing
a nondegenerate integral element (in $\Hom ({}_BA_B, {}_BB_B)$
such as a Frobenius homomorphism).  However, if $R = Z(A)$, as mentioned
above an antipode $\tau_T$ is definable in the same way as $\tau = \tau_S$.  

The computations in this section expose the hypotheses that 
are necessary for the antipode defined in \cite[4.4]{KN}.  The correspondence
of the theory in this section with that in \cite{KN} is discussed
in \cite[8.9]{KS}, and depends on the equation for a depth
two Frobenius extension $A | B$ with trivial one-dimensional centralizer:  
\begin{equation}
[t, \alpha ] = E_M E_{M_1}(\psi(t) e_1 e_2 \phi(\alpha))
\end{equation}
for certain anti-isomorphisms $\psi: T \rightarrow C_{M_2}(M)$
defined  in \cite[8.2]{KS} and $\phi: S \rightarrow C_{M_1}(N)$ defined
in \cite[8.4]{KS}.  

\section{The $\pi$-method for depth two
extensions}

In this section we extend Doi and Takeuchi's $\pi$-method for Hopf-Galois extensions \cite{DT} to D2 extensions.
Then we extend the antipode in the previous section
to a certain bimodule hom-group for irreducible extensions. We point
out that the $\pi$-method yields a nullhomotopic mapping between relative
Hochschild cochains with coefficients
for the irreducible extension $A \| B$ and a certain Hochschild 
cohomology theory with coefficients for the $R$-coring $T$.   

Suppose $A \| B$ is an rD2 extension,
and ${}_AM$ is a module.  Again let
$T = (A \o_B A)^B$, the right bialgebroid
over the centralizer $R = A^B$.  Let
$u_j \in T$ and $\gamma_j \in S$ be rD2
quasibases.  Recall from \cite{LK2006A}
the coaction $\rho^T: A \to A \o_R T$ which makes
$A$ into a right comodule algebra: in Sweedler notation this is given by 
\begin{equation}
a\0 \o_R a\1 = \sum_j \gamma_j(a) \o_R u_j
\end{equation}
Clearly, $b\0 \o b\1 = b \o 1_T$ if
$b \in B$. 
 
\begin{prop}
The mapping $\pi_L: \Hom ({}_RT, {}_RM) \longrightarrow \Hom ({}_BA, {}_BM)$
given by
\begin{equation}
\pi_L(f)(a) = a\0 f(a\1)
\end{equation}
is a left $B$-linear, $R$-linear isomorphism.  
\end{prop}
\begin{proof}
The inverse to $\pi_L$ is given by
\begin{equation}
\label{eq: pi inverse}
\pi_L^{-1}(g)(t) = t^1 g(t^2)
\end{equation}
for $g \in \Hom ({}_BA, {}_BM)$.  We note
that 
$$ t^1 \pi_L(f)(t^2) = \sum_j t^1 \gamma_j(t^2) f(u_j) = f(t)$$
for $f \in \Hom ({}_RT, {}_RM)$,
$t \in T$, since $t^1 \gamma_j(t^2) \in A^B = R$ and $\sum_j t^1 \gamma_j(t^2)
u_j = t$.  We also note that
$$ \sum_j \gamma_j(a) u^1_j g(u^2_j) = g(a)
$$
for $a \in A, g \in \Hom ({}_BA, {}_BM)$
by eq.~(\ref{eq: right d2 qb}). 
Hence, $\pi_L^{-1}$ is indeed the inverse
of $\pi_L$.  

The mapping $\pi_L$ is left $B$-linear,
since 
$$\pi_L(bf)(a) = a\0 b f(a\1) = 
\pi(f) (ab).$$ 
Note that $\pi_L^{-1}$ is left $R$-linear
since $$\pi_L^{-1}(rg)(t) = t^1 (rg)(t^2)
= t^1 g(t^2r) = \pi_L^{-1}(g)(tr). \qed $$ 
\renewcommand{\qed}{}\end{proof}

Similarly, if $A \| B$ is $\ell$D2,
$N_A$ is a module and $\beta_i \in S$,
$t_i \in T$ are $\ell$D2 quasibases,
then we have the right module dual of
the proposition:  
\begin{equation}
\pi_{\mathcal{R}}: \ \Hom (T_R, N_R) \stackrel{\cong}{\longrightarrow}
\Hom (A_B, N_B), \ \  \pi_{\mathcal{R}}(h)(a) = \sum_i h(t_i) \beta_i(a),
\end{equation}
which has inverse mapping given for $g \in \Hom (A_B, N_B), t\in T$ by
\begin{equation}
\pi_{\mathcal{R}}^{-1}(g)(t) = g(t^1)t^2
\end{equation}

Suppose that ${}_AP_A$ is a bimodule.    We let $P^B$ denote the elements $p \in P$
satisfying $pb = bp$ for all $b \in B$. Note that
$P^B$ is a natural $(R,R)$-bimodule, since $R$ and $B$ commute.

First note that $\pi_L: \Hom ({}_RT, {}_RP) \stackrel{\cong}{\longrightarrow}
\Hom ({}_BA, {}_BP)$ restricts to
\begin{equation}
\pi_L:\ \Hom ({}_RT, {}_RP^B) \stackrel{\cong}{\longrightarrow}
\Hom ({}_BA_B, {}_BP_B) 
\end{equation}
Similarly $\pi_{\mathcal{R}}$ above restricts to
\begin{equation}
\pi_{\mathcal{R}}: \Hom (T_R, P^B_R) \stackrel{\cong}{\longrightarrow} \Hom ({}_BA_B, {}_BP_B).
\end{equation}
Note that the inverses are given by
\begin{equation}
\pi_{\mathcal{R}}^{-1}(h)(u) = h(u^1)u^2, \ \ \pi_L^{-1}(h)(u) = u^1 h(u^2)
\end{equation}
for $h \in \Hom ({}_BA_B, {}_BP_B), u \in T$.  Of course if $P = A$ these restricted
mappings recover the isomorphisms
$\Hom (T_R, R_R) \cong S \cong \Hom ({}_RT, {}_RR)$ below eq.~(\ref{eq: bra-ket pairing}).

Now suppose that the ring extension $A \| B$ is irreducible.
The antipode in
the previous sections then extends to a bijection of the hom-group $\Hom ({}_BA_B, {}_BP_B)$ onto itself as follows.
Since $R \subseteq B$, it follows that $\Hom (T_R, P^B_R) = \Hom ({}_RT, {}_RP^B)$.  Then define the mapping 
$$ \sigma_P = \pi_L \circ \pi_{\mathcal{R}}^{-1}: \Hom ({}_BA_B, {}_BP_B) \stackrel{\cong}{\longrightarrow} \Hom ({}_BA_B, {}_BP_B) $$
which is then given by 
\begin{equation}
\sigma_P(\alpha)(a) = \pi_L( \pi_{\mathcal{R}}^{-1}(\alpha))(a) = \sum_j \gamma_j(a)
\alpha(u_j^1)u_j^2 
\end{equation}
for $\alpha \in \Hom ({}_BA_B, {}_BP_B)$.  Now let $\alpha^{\sigma} = 
\sigma_P(\alpha)$.  Then for all $t \in T$, 
\begin{equation}
\label{eq: antipode extended}
\alpha(t^1)t^2 = t^1 \alpha^{\sigma}(t^2),
\end{equation}
which follows from the following short computation:
$$ t^1 \sigma_P(\alpha)(t^2) =  \sum_j t^1\gamma_j(t^2)
\alpha(u_j^1)u_j^2 = \alpha(t^1)t^2,$$
since $t^1 \gamma_j(t^2) \in R \subseteq B$ and eq.~(\ref{eq: right d2 qb}). 
A glance at eq.~(\ref{eq: tau alt}) shows that $\sigma_A = \tau$,
the antipode of $S$ defined in previous sections. 
The proof of the proposition is similar to  previously 
and therefore omitted.

\begin{prop}
Suppose $A \| B$ is an irreducible extension and ${}_AP_A$ is a bimodule.
Define two pairings of $\Hom ({}_BA_B, {}_BP_B)$ and $T$ with
values in $P^B$ by $\bra \alpha, t \ket = \alpha(t^1)t^2$
and $[u, \beta ] = u^1 \beta(u^2)$.  Then the two pairings are non-degenerate
w.r.t.\ $\alpha, \beta \in \Hom ({}_BA_B, {}_BP_B)$.  Moreover, the mapping
 $\alpha \mapsto \alpha^{\sigma}$  defined by
$\bra \alpha, t \ket = [t, \alpha^{\sigma}]$ is a bijection
satisfying $\sigma^{\pm 1}(\lambda_p) = \rho_p$ for
each $p \in P$, 
and $\alpha^{\sigma}(1_A) = \alpha(1_A)$.  
\end{prop} 

\subsection{Remark on  relative Hochschild cohomology with coefficients}
We remark below on how the $\pi$-method leads to a nullhomotopic
mapping between certain Hochschild cohomology theories.  
Continuing the notation just above, note
that $\Hom ({}_BA_B, {}_BP_B)$ is the
first relative Hochschild cochain group
of $A \| B$ with coefficients in a bimodule ${}_AP_A$, denoted by $C^1(A,B;P)$.  The zero'th group
is $P^B$ with differential $d^0: P^B \to
\Hom ({}_BA_B, {}_BP_B)$ given
by $d^0(p) = \rho_p - \lambda_p$
for $p \in P^B$.  

The second relative Hochschild cochain group
is $C^2(A, B; P) = \Hom ({}_BA \o_B A_B,
{}_BP_B)$. The differential at this level
is given by $d^1: \Hom ({}_BA_B, {}_BP_B)
\to \Hom ({}_BA\o_B A_B, {}_BP_B)$
defined by
\begin{equation}
(d^1f)(x \o_B y) = xf(y) - f(xy) + f(x)y
\end{equation}
The cohomology groups are denoted by $HH^n(A,B;P)$ for $n \geq 0$; recall or 
note that $HH^1(A,B;P)$ is isomorphic to the group of derivations
killing $B$ 
modulo the group of inner derivations w.r.t.\ elements in $P^B$.    
For the sake of brevity we refer
the reader to textbooks on homological algebra for the details of higher order cochain groups and differentials.   

Note that applications of
the hom-tensor relation to the left- and right-handed
$\pi$-methods with ${}_AM_A = 
\Hom (A_B, P_B)$ yields two isomorphisms
\begin{equation}
\Hom ({}_BA \o_B A_B, {}_BP_B) \stackrel{\cong}{\longrightarrow} \Hom ({}_RT \o_R T, {}_RP^B_R)
\end{equation}
denoted by $\pi^{2}_L$ and $\pi^{2}_{\mathcal{R}}$ given
by 
\begin{eqnarray}
\pi^{2}_L(h)(u \o_R t) &=& u^1 h(u^2t^1 \o_B t^2) \\
\pi^{2}_{\mathcal{R}}(h)(u \o_R t) & = & h(u^1 \o_B u^2 t^1) t^2.
\end{eqnarray}
The inverses are given by 
\begin{eqnarray}
\pi^{-2}_L(g)(x \o_B y) & = & \sum_{j,k} \gamma_k(x \gamma_j(y))g(u_k \o_R u_j) \\
\pi^{-2}_{\mathcal{R}}(g)(x \o_B y) & = & \sum_{i,j} g(t_i \o_R t_j)\beta_j(\beta_i(x)y) 
\end{eqnarray}
Likewise we define the obvious generalized isomorphisms
$\pi^n_L$ and $\pi^n_{\mathcal{R}}$ on
the $n$-cochains $C^n(A,B;P)$.  

Brzezinski and Wisbauer define a Hochschild
cohomology of an $R$-coring with coefficients in an $(R,R)$-bimodule \cite[30.15]{BW}.  For the $R$-coring $T$
and $(R,R)$-bimodule $P^B$, this specializes to the zero'th cochain group
 $\Hom ({}_RT_R,
{}_RP^B_R)$ (which is equal to both one-sided $R$-linear hom-groups considered
above since $A \| B$ is irreducible),
and first 
cochain group $\Hom ({}_RT \o_R T_R, {}_RP^B_R)$. Let $\eps: T \to R$ be the counit of $T$ given by $\eps(t) = t^1t^2$,
the multiplication mapping $A \o_B A \to A$
restricted to $T$.  The diffential is given
by ($h \in \Hom ({}_R T_R, {}_RP^B_R)$)
\begin{equation}
(\delta^0 h)(u \o_R t) = \eps(u) h(t) - h(u) \eps(t) 
\end{equation}
\begin{equation}
(\delta^1 g)(u \o_R v \o_R t) =
\eps(u) g(v \o_R t) - g(u \eps(v) \o_R t)
+ g(u \o_R v) \eps(t) 
\end{equation}
for $g \in \Hom ({}_RT \o_R T_R, {}_RP^B_R)$.
For example, $\delta^1 \delta^0 h = 0$
by a short computation.
The higher cochain groups and differentials
are defined similarly and we refer to \cite[30.15]{BW} for the details. (This cochain
complex recovers Hochschild relative
cochains in case the $A$-coring $A \o_B A$
with counit $\eps'(x \o y) = xy$
takes the place of $R$, $T$ and $\eps$.)
Denote the cochain groups in this complex
by $C^n(T,R; P^B)$, and call it the
Hochschild coring complex.  

Next we define a cochain homomorphism
$\Phi_n: C^n(A,B;P) \to C^{n-1}(T,R;P^B)$,
for $n \geq 1$, 
noting the shift of one downwards
 in degree in the Hochschild coring complex. Define $\Phi_1: \Hom ({}_BA_B, {}_BP_B) \to \Hom ({}_RT_R,{}_RP^B_R)$
as 
\begin{equation}
\Phi_1 = \pi_L^{-1}  + \pi_{\mathcal{R}}^{-1}, \ \ \Phi_1(h)(u) = u^1h(u^2) + h(u^1)u^2 
\end{equation}
for $h \in \Hom ({}_BA_B, {}_BP_B)$. Note that $\Phi_1$ kills $B$-linear
derivations. 

Define $\Phi_2: \Hom ({}_BA \o_B A_B, {}_BP_B) \to \Hom ({}_R T \o_R T_R, {}_RP^B_R)$ by $\Phi_2 = \pi^2_L - \pi^2_{\mathcal{R}}$, or  in more detail, 
\begin{equation}
(\Phi_2 g)(u \o_R t) = u^1 g(u^2 t^1 \o_B t^2) - g(u^1 \o_B u^2 t^1)t^2
\end{equation}
The $n$'th mapping $\Phi_n$ is easily
defined from two obvious generalized  mappings, $\pi^n_L, \pi_{\mathcal{R}}^n: \Hom ({}_BA \o_B \cdots \o_B A_B, {}_BP_B) \to \Hom ({}_RT \o_R \cdots \o_R T_R, {}_RP^B_R)$; as $\Phi_n = \pi^n_L + (-1)^{n+1}\pi^n_R$.  

We compute for $f \in \Hom ({}_BA_B, {}_BP_B)$,
 $$\delta^0 \Phi_1 f (u \o_R t) =
(\Phi_1 f)(u^1u^2 t) - (\Phi_1 f)(u t^1 t^2) $$
$$ = u^1u^2 t^1 f(t^2) + f(u^1u^2 t^1)t^2
-f(u^1)u^2 t^1 t^2 - u^1 f(u^2 t^1 t^2) $$
$$ = u^1 (df)(u^2 t^1 \o_B t^2) - (df)(u^1 \o_B u^2 t^1)t^2 = (\Phi_2 d^1 f)(u \o_R t)$$
after two middle terms cancel. 

Moreover, for $g \in \Hom ({}_BA \o_B A_B, {}_BP_B)$,
$$ \Phi_3(dg)(v \o_R u \o_R t) = v^1 (dg)(v^2u^1 \o_B u^2 t^1 \o_B t^2) + (dg)(v^1 \o_B v^2 u^1 \o_B u^2 t^1) t^2 $$
$$ = (\Phi_2 g)(\eps(v)u \o_R t) - (\Phi_2 g)(v \eps(u) \o_R t) + (\Phi_2 g)( v \o_R u \eps(t)) = (\delta^1 \Phi_2 g)(v \o_R u \o_R t)$$
after cancellation of the pair of middle  terms $\pm v^1 g(v^2u^1 \o_B u^2 t^1)t^2$.  

We omit the tedious but similar computation 
in degree $n$ which establishes that $\Phi$
is a cochain mapping. 
\begin{prop}
Suppose $A \| B$ is an irreducible extension
and ${}_AP_A$ is a bimodule.
Let $R = A^B$ and $T = (A \o_B A)^B$. Then the mapping of cochain groups $\Phi_n: C^n(A,B;P) \rightarrow C^{n-1}(T,R;P^B)$
is nullhomotopic.    
\end{prop}
\begin{proof}
We define a homotopy $s_n: C^{n+2}(A,B;P) \to C^{n}(T,R:P^B)$ first in degree zero by $s_0(f)(t) = f(t^1 \o_B t^2)$ where $f \in \Hom ({}_BA \o_B A_B, {}_BP_B)$, $t \in T \subseteq A \o_B A$, so $f(t) \in P^B$. We claim there is a natural
inclusion $\iota_n$ of $T \o_R \cdots \o_R T$
($n$ times $T$) into $A \o_B \cdots \o_B A$
($n+1$ times $A$) given by
\begin{equation}
\iota_n(u_1 \o_R \cdots \o_R u_n) = 
u^1_1 \o_B u^2_1 u^1_2 \o_B \cdots
\o_B u^2_{n-1} u_n^1 \o_B u_n^2 
\end{equation}
 In fact, $\iota_n$ is an isomorphism onto
$(A \o_B \cdots \o_B A)^B$ which follows
from showing that for any ring $C$ and  $(A,C)$-bimodule $M$
\begin{equation}
T \o_R M \stackrel{\cong}{\longrightarrow}
A \o_B M
\end{equation}
 via $t \o_R m \mapsto t^1 \o_B t^2m$.  
This is a $(B,C)$-bimodule isomorphism with inverse given by $a \o_B m \mapsto
\sum_i t_i \o_R \beta_i(a)m$ using left
D2 quasibases $\beta_i \in S$, $t_i \in T$.
Now we may apply this with $C = A$, $M =  A \o_B A, A \o_B A \o_B A, \ldots$, then restrict to $(-)^B$, substitute
and iterate to prove the claim. 

Now define $s_n(g) = g \circ \iota_{n+1}$
for $g \in C^{n+2}(A,B;P)$.  
We note that $\delta^n s_n + s_{n+1} d^{n+2} = \Phi_{n+2}$;
e.g., for $f \in C^2(A,B;P)$, $u,t \in T$, 
$$\delta^0 (s_0 f)(u \o_R t) + s_1(d^2f)(u \o_R t) = \eps(u) f(t) 
- f(u) \eps(t) + (df)(u^1 \o_B u^2 t^1 \o_B t^2) $$
$$ = u^1 f(u^2 t^1 \o_B t^2)  - f(u^1 \o_B u^2t^1)t^2  = (\Phi_2 f)(u \o_R t)
$$
after cancellation of the middle terms in
$(df)(u^1 \o_B u^2 t^1 \o_B t^2) $. 
The general case is computed similarly.
\end{proof}


\section{Appendix: axioms for Hopf algebroids}

In this appendix we review the definitions of left bialgebroid, Lu's Hopf
algebroid, right bialgebroid and the definition of Hopf algebroid by
B\"ohm-Szlach\'anyi.  
First, for the definition of a left bialgebroid $(H,R,s_L,t_L,\cop, \eps)$,  $H$ and $R$ are $K$-algebras and all 
maps are $K$-linear.  First,
recall from \cite{Lu}  that
the \textit{source} and \textit{target} maps $s_L$ and $t_L$ 
are algebra homomorphism and anti-homomorphism, respectively,
of $R$ into $H$ such that $s_L(r)t_L(s) = t_L(s)s_L(r)$ 
for all $r,s \in R$.
This induces an $R$-$R$-bimodule structure on $H$ (from the left in this case) by
$r \cdot h \cdot s = s_L(r)t_L(s) h$ ($h \in H$).  With respect to this bimodule
structure, $(H,\cop,\eps)$ is an $R$-coring (cf.\ \cite{Sw}), i.e. with coassociative
coproduct and $R$-$R$-bimodule map $\cop: H \to H \o_R H$ and
counit $\eps: H \to R$ (also an $R$-bimodule mapping).  The image
of $\cop$, written in Sweedler notation, is required to satisfy
\begin{equation}
 a\1 t_L(r) \o a\2 = a\1 \o a\2 s_L(r)
\end{equation}
for all $a \in H, r \in R$.  It then makes sense to require that $\cop$ be homomorphic:
\begin{equation}
\cop(ab) = \cop(a) \cop(b), \ \ \ \cop(1) = 1 \o 1
\end{equation}
for all $a,b \in H$. The counit must satisfy the following modified augmentation law:
\begin{equation}
\eps(ab) = \eps(a s(\eps(b))) = \eps(a t(\eps(b))), \ \ \eps(1_H) = 1_R. 
\end{equation}

The axioms of a right bialgebroid $H'$ are opposite those of a left bialgebroid in
the sense that $H'$ obtains its $R$-bimodule structure from the right via its
source and target maps and, from the left bialgebroid $H$ above, we have that 
 $(H^{\rm op}, R, t_L^{\rm op}, s_L^{\rm op}, \cop, \eps)$
(in that precise order) is a right bialgebroid: 
for the explicit axioms, see \cite[Section~2]{KS}. 

In addition, the left $R$-bialgebroid $H$
is a \textit{Hopf algebroid} in the sense of Lu   $(H,R,\tau)$ 
if (antipode) $\tau: H \to H$ is an algebra anti-automorphism such that 
\begin{enumerate}
\item $\tau t_L = s_L$;
\item $\tau(a\1)a\2 = t_L (\eps (\tau(a)))$ for every $a \in A$;
\item there is a linear section $\eta: H \o_R H \to H \o_K H$ to the natural
projection $H \o_K H \to H \o_R H$ such that:
$$ \mu (H \o \tau) \eta \cop = s_L \eps. $$
\end{enumerate}

The following is one of several equivalent definitions of
B\"ohm-Szlach\'anyi's Hopf algebroid \cite{BS}, excerpted from \cite[8.7]{KS}.
\begin{defi}
We call $H$ a Hopf algebroid if there are left and right bialgebroid structures $(H,R, s_L, t_L, \cop_L, \eps_L)$ and $(H,R^{\rm op}, s_R, t_R, \cop_R, \eps_R)$ such that
\begin{enumerate}
\item $\Im s_R = \Im t_L$ and $\Im t_R = \Im s_L$,
\item $(1 \o \cop_L)\cop_R = (\cop_R \o 1) \cop_L$ and $( 1 \o \cop_R)\cop_L = (\cop_L \o 1)\cop_R$
\end{enumerate}
with anti-automorphism $\tau: H \to H$ (called an antipode) such that
\begin{description}
\item[(i)] $\tau(at_R(r)) = s_R(r) \tau(a)$ and $\tau(t_L(r)a) = \tau(a)s_L(r)$ for $r \in R, a \in H$,
\item[(ii)] $a\i \tau(a\j) = s_L(\eps_L(a))$ and $\tau(a\1)a\2 = s_R(\eps_R(a))$
\end{description}
where $\cop_R(a) = a\i \o a\j$ and $\cop_L(a) = a\1 \o a\2$. 
\end{defi}

The relationship between Lu's Hopf algebroid and this alternative Hopf algebroid with more pleasant tensor  categorical properties is discussed
in \cite{BS, B} and other papers by B\"ohm and Szlach\'anyi.

\end{document}